\documentclass{article}[13pt]

\usepackage{graphicx}
\usepackage{amsthm,amsmath}

\numberwithin{equation}{section}
\theoremstyle{plain}
\newtheorem{thm}{Theorem}

\newtheorem{lem}{Lemma}

\newcommand{\wht}{\widehat}
\newcommand{\wtd}{\widetilde}

\newcommand{\wbox}{\sqcap\llap{$\sqcup$}}

\newcommand{\bmu}{\boldsymbol{\mu}}
\newcommand{\btheta}{\boldsymbol{\theta}}
\newcommand{\blambda}{\boldsymbol{\lambda}}

\begin{document}

\begin{center}
{\large \bf OPTIMAL UCB ADJUSTMENTS FOR LARGE ARM SIZES}


\medskip
{\large by Hock Peng Chan and Shouri Hu}

\medskip
{\large National University of Singapore}



\bigskip
{\bf Abstract}
\end{center}

\smallskip
The regret lower bound of Lai and Robbins (1985), 
the gold standard for checking optimality of bandit algorithms,
considers arm size fixed as sample size goes to infinity.
We show that when arm size increases polynomially with sample size,
a surprisingly smaller lower bound is achievable. 
This is because the larger experimentation costs when there are more arms permit regret savings by exploiting the best performer more often.
In particular we are able to construct a UCB-Large algorithm that adaptively exploits more when there are more arms.
It achieves the smaller lower bound and is thus optimal.
Numerical experiments show that UCB-Large performs better than classical UCB that does not correct for arm size,
and better than Thompson sampling.




\section{Introduction}

Let there be $K$ arms (populations) from which rewards (observations) are drawn. 
The multi-armed bandit problem is the design of sequential samplers that allocate sampling to maximize expected sum of rewards.

Consider a family of densities $\{ f_{\theta}: \theta \in \Theta \}$,
with respect to a measure on the real line.
Let rewards from arm $k$ be independent and identically distributed (i.i.d.) with density $f_{\theta_k}$,
for unknown $\theta_k$.
Let $P_{\btheta}$ ($E_{\btheta}$) denote probability (expectation) with respect to $\btheta (=\btheta_K) = (\theta_1, \ldots, \theta_K)$.
Let $\mu(\theta)$ be the mean of $f_{\theta}$,
and let $\mu_* = \max(\mu_1,\ldots,\mu_K)$,
where $\mu_k = \mu(\theta_k)$.

Maximizing expected sum of $N$ rewards is the same as minimizing the regret
$$R_N(\btheta) := \sum_{k=1}^K (\mu_*-\mu_k) E_{\btheta} N_k, 
$$
where $N_k$ is the number of rewards from arm $k$.
Let
\begin{equation} \label{Lth}
r(\btheta) = \sum_{k: \mu_k < \mu_*} \tfrac{\mu_*-\mu_k}{D(\theta_k|\theta_*)},
\end{equation}
where $D(\theta|\lambda) = E_{\theta} [\log \tfrac{f_{\theta}(X)}{f_{\lambda}(X)}]$ is the Kullback-Leibler (KL)-information number between $f_{\theta}$ and $f_{\lambda}$,
and $\theta_*$ is such that $\mu(\theta_*) = \mu_*$.
The celebrated lower bound result of Lai and Robbins (1985) is that if the regret $R_N(\btheta)$ grows at a sub-polynomial rate with respect to $N$ for each $\btheta$ (with $K$ fixed),  
then as $N \rightarrow \infty$,
\begin{equation} \label{RNt}
R_N(\btheta) \geq [1+o(1)] r(\btheta) \log N. 
\end{equation}

Lai (1987), 
Agrawal (1995), 
Burnetas and Katehakis (1996) and Capp\'{e} et al. (2013) constructed upper confidence bound (UCB) algorithms that have regret achieving equality in (\ref{RNt}) on exponential families,
and are thus optimal.
UCB-Agrawal (Burnetas and Katehakis suggested the same algorithm) improves upon UCB-Lai in not requiring advance knowledge of $N$.
Auer,
Cesa-Bianchi and Fischer (2012) provided finite $N$ upper bounds of UCB for bounded rewards.
Chan (2019) showed that instead of applying confidence bounds that are specific to a given exponential family,
subsampling can be applied to achieve optimality on {\it unspecified} exponential families.

We show here that if $K=N^{\zeta+o(1)}$ for some $0 \leq \zeta < 1$ as $N \rightarrow \infty$,
then instead of (\ref{RNt}) we have
\begin{equation} \label{R2}
R_N(\btheta_K) \geq [1-\zeta+o(1)] r(\btheta_K) \log N.
\end{equation}
The smaller lower bound when $\zeta>0$ is not due to technical difficulties in extending the lower bound proof of Lai and Robbins.
Rather we show that it is sharp, 
by constructing a UCB-Large bandit algorithm that achieves this smaller lower bound.
In addition we are able to overcome the technical difficulties mentioned in Burnetas and Katehakis to show that UCB is optimal when rewards are normal with unknown and unequal variances. 

The improvements of UCB-Large over classical UCB strengthen UCB as a competitor of Bayesian approaches to the multi-armed bandit problem,
see Gittins (1979),
Gittin and Jones (1979),
Brezzi and Lai (2000) as well as Thompson (1933),
Berry and Fristedt (1985),
Kaufmann, Capp\'{e} and Munos (2012) and Korda, Kaufmann and Munos (2012).
The improvements are due to adaptations of UCB to take into account the unavoidable experimentation costs unique to a particular problem,
in this case the higher costs when the number of arms is large.
The construction of optimal bandit algorithms for irreversible rules in Hu and Wei (1989) is also based on this principle.

Algorithms for large arm sizes have been constructed in Berry et al. (1995) and Bonald and Prouti\`{e}re (2013) for Bernoulli rewards,
and in Chan and Hu (2019) for general rewards that are bounded above.
A key difference of these algorithms is that they assume infinite number of arms are available, 
so that not all arms can be sampled.

The layout of this paper is as follows.
In Section 2 we propose UCB-Large for general exponential families.
In Section 3 we show optimality for normal rewards,
for both variances known and unknown.
The restriction to normal rewards is to avoid technical complexities that occur with unbounded arm means as the number of arms goes to infinity,
and with the gap between the optimal and best inferior arm going to zero. 
In Section 4 we consider gaps that decrease polynomially fast.
In Section 5 we confirm, 
via numerical studies, 
the improvements of UCB-Large over classical UCB.
In Sections 6--8 we prove the results of Sections 3 and 4.

\section{UCB-Large}

Consider the one-dimensional exponential family
\begin{equation} \label{expo}
f_{\theta}(x) = e^{\theta x-\psi(\theta)} f_0(x), \quad \theta \in \Theta,
\end{equation}
where $\psi(\theta) = \log E_0 e^{\theta X}$ and $\Theta = \{ \theta: \psi(\theta) < \infty \}$.
Let $\theta_x$ be such that $\psi'(\theta_x)=x$.
Under (\ref{expo}), 
the large deviations rate function
$$I_u(x) = (\theta_x-\theta_u) x-[\psi(\theta_x)-\psi(\theta_u)] = D(\theta_x|\theta_u).
$$

Let $S_t = \sum_{s=1}^t X_s$,
$S_{kt} = \sum_{s=1}^t X_{ks}$,
$\bar X_t = \frac{S_t}{t}$ and $\bar X_{kt} = \frac{S_{kt}}{t}$.
Let $U_t(\bar X_t,b)$ be the upper confidence bound of $(X_1,\ldots,X_t)$,
with respect to confidence coefficient $b$,
where 
\begin{equation} \label{ruX}
U_t(x,b) = \inf \{ u \geq x: t I_u(x) \geq b \}.
\end{equation}
Let $b_n$ be non-negative and monotone increasing for $n \in [1,\infty)$.
Let $n_k$ be the number of rewards from arm $k$ when there are $n$ total rewards.
Let $n_1=\cdots=n_K=1$ when $n=K$,
that is we initialize with one reward allocated to each arm.
Agrawal (1995) proposed the following UCB procedure.

\smallskip
\underline{UCB-Agrawal}.
For $n=K, \ldots,N-1$:
Compute the confidence bounds
\begin{equation} \label{Uk}
U_k^n = U_{n_k}(\bar X_{kn_k},b_n), \quad 1 \leq k \leq K,
\end{equation}
and sample from the arm with largest confidence bound.

\smallskip
Agrawal showed that UCB-Agrawal achieves the Lai-Robbins lower bound (\ref{RNt}) (with $K$ fixed), 
when $b_n \sim \log n$ with
\begin{equation} \label{bn}
b_n - \log n- \log \log n \rightarrow \infty \mbox{ as } n \rightarrow \infty.
\end{equation}
Burnetas and Katehakis (1996) showed that UCB-Agrawal achieves (\ref{RNt}) for $b_n=\log n$, 
under a condition that they remarked is satisfied for normal densities with known variances,
and claimed that their analysis carries over to general $b_n \sim \log n$.
Motivated by the multi-parameter regret lower bounds in Burnetas and Katehakis (1996),
Honda and Takemura (2010) constructed an asymptotically optimal DMED algorithm for distributions with bounded support.
Capp\'{e} et al.(2013) provided finite $N$ regret for their KL-UCB algorithm when $b_n = \log n + 3 \log \log n$ but recommended $b_n = \log n$ for practical use.

In practice $K$ can be large.
We show in Section 3 that optimality is extended to $K \rightarrow \infty$ by simply replacing $b_n$ in (\ref{Uk}) with $b_{n/K}$.

\smallskip
\underline{UCB-Large}.
For $n=K, \ldots, N-1$:
Compute the confidence bounds
\begin{equation} \label{Ukn}
U_k^n = U_{n_k}(\bar X_{kn_k},b_{n/K}), \quad 1 \leq k \leq K, 
\end{equation}
and sample from the arm with largest confidence bound.

\smallskip
In addition to the most natural $b_m=\log m$,
our numerical studies in Section 5 include $b_m = \chi \log m$ for $\frac{1}{2} \leq \chi < 1$ and $b_m = \log(e-1+m)-\sqrt{\log(e-1+m)}$.
These confidence coefficients are justified in Sections 3 and 4.
The examples below cover the most important exponential families.

\smallskip
{\sc Example} 1. 
Consider $f_{\mu}$ the normal density with mean $\mu$ and variance~1.
Here $I_u(x) = \frac{(u-x)^2}{2}$.
The confidence bounds under UCB-Agrawal for $b_n=\log n$ are
\begin{equation} \label{Agr}
U_k^n = \bar X_{kn_k} + \sqrt{\tfrac{2 \log n}{n_k}}.
\end{equation}
The confidence bounds under UCB-Large for the same $b_n$ is 
\begin{equation} \label{L1}
U_k^n = \bar X_{kn_k} + \sqrt{\tfrac{2 \log(n/K)}{n_k}}.
\end{equation}

The intuition behind (\ref{L1}) is as follows.
The confidence bounds (\ref{Agr}) are designed so that the exploitation cost is $o(\log N)$.
The exploitation cost is the cost of sampling the arm with largest sample mean when it is in fact an inferior arm.
For large $K$ this control is overly strict as the exploration cost, 
of order $K \log N$,
is much larger and so for optimality the UCB should reduce exploration up to the point where exploitation cost reaches $o(K \log N)$.
This is achieved by the insertions of $K$ in the confidence bounds (\ref{L1}).

\smallskip
{\sc Example} 2. 
Consider normal rewards with unknown and unequal variances.
Here UCB is extended to a two-dimensional exponential family.
Let $\theta=(\mu,\sigma^2)$ and $\Theta = \{ \theta: \sigma^2>0 \}$.
For $\theta \in \Theta$,
let
\begin{equation} \label{tdisn}
f_{\theta}(x) = \tfrac{1}{\sigma \sqrt{2 \pi}} e^{-\tfrac{(x-\mu)^2}{2 \sigma^2}}. 
\end{equation}
Let $M(z) = \frac{1}{2} \log(1+z^2)$.

Burnetas and Katehakis (1996) proposed upper confidence bounds
\begin{eqnarray} \label{U4}
U_k^n & = & \inf \{ u \geq \bar X_{kn_k}: n_k M(\tfrac{u-\bar X_{kn_k}}{\hat \sigma_{kn_k}}) \geq \log n \} \\ \nonumber
& = & \bar X_{kn_k} + \wht \sigma_{kn_k} \sqrt{\exp(\tfrac{2 \log n}{n_k})-1},
\end{eqnarray}
where $\wht \sigma^2_{ks} = s^{-1} \sum_{t=1}^s (X_{kt}-\bar X_{ks})^2$.

They showed that (for $K$ fixed) if the regret grows sub-polynomially with $N$ for each $\btheta$,
then as $N \rightarrow \infty$,
\begin{equation} \label{lowert}
R_N(\btheta) \geq [1+o(1)] r(\btheta) \log N, \mbox{ where } r(\btheta) = \sum_{k: \mu_k < \mu_*} \tfrac{\mu_*-\mu_k}{M(\frac{\mu_*-\mu_k}{\sigma_k})}.
\end{equation}
They did not show that (\ref{U4}) has regret achieving the lower bound in (\ref{lowert}), 
due to difficulties with the tail probabilities of non-central $t$-distributions.
We overcome these difficulties (and extend to $K$ large) by applying instead the confidence bounds
\begin{eqnarray} \label{UCBt}
U_k^n & = & \inf\{ u \geq \bar X_{kn_k}: (n_k-1) M(\tfrac{u-\bar X_{kn_k}}{\hat \sigma_{kn_k}}) \geq b_{n/K} \} \\ \nonumber
& = & \bar X_{kn_k} + \wht \sigma_{kn_k} \sqrt{\exp(\tfrac{2b_{n/K}}{n_k-1})-1}.
\end{eqnarray}
The subtraction of 1 from $n_k$ in (\ref{UCBt}) can be viewed as the effective sample size reduction to account for the estimation of $\sigma_k^2$.



\smallskip
{\sc Example} 3.
Consider $f_{\theta}$ the Bernoulli density (with respect to counting measure on $\{ 0,1 \}$) with mean $\mu(\theta)$. 
The large deviations rate function
\begin{equation} \label{LDB}
I_u(x)= x \log(\tfrac{x}{u}) + (1-x) \log (\tfrac{1-x}{1-u}), \quad 0 \leq x \leq 1,
\end{equation}
with $0 \log 0=0$.
The confidence bound $U_k^n$ is the larger root in $u$ of $I_u(x_k)=y_k$,
with $x_k = \bar X_{kn_k}$ and $y_k=\frac{b_{n/K}}{n_k}$. 
A quick way to compute $U_k^n$ is to initialize with $v_{k0} \in [x_k,1]$ and solve iteratively, 
for $i \geq 0$,
\begin{equation} \label{xv}
x_k \log(\tfrac{x_k}{v_{ki}})+(1-x_k) \log(\tfrac{1-x_k}{1-v_{k,i+1}})=y_k.
\end{equation}
A computational advantage of (\ref{xv}) is that the iterations
$$v_{k,i+1} = 1-(d_k/v_{ki}^{x_k})^{\frac{1}{1-x_k}},
$$
with $d_k = x_k^{x_k} (1-x_k)^{1-x_k} e^{-y_k}$,
can be executed simultaneously on all arms, 
by common operations on $(v_{1i}, \ldots, v_{Ki})$.

\section{Regret lower bound and optimality of UCB-Large for normal rewards}

Let $a^+ = \max(a,0)$ and let $J(\bmu) = \# \{ k: \mu_k < \mu_* \}$ be the number of inferior arms with respect to $\bmu$.
We say that $\Delta_K \rightarrow 0$ at a sub-polynomial rate if $\Delta_K K^{\epsilon} \rightarrow \infty$ for all $\epsilon > 0$.
We consider either $K=N^{\zeta+o(1)}(\rightarrow \infty)$ for some $0 < \zeta < 1$ or $K$ fixed (i.e. $\zeta=0$) as $N \rightarrow \infty$.

\subsection{Normal rewards with unit variances}

Let $X_{kt} \stackrel{\rm i.i.d.}{\sim}$ N($\mu_k,1)$,
$t \geq 1$,
be the rewards of arm $k$.
Let
$$\Theta(\Delta_K) = \{ \bmu_K: \max_{k: \mu_k < \mu_*} \mu_k \leq \mu_*-\Delta_K, \ \max_k |\mu_k| \leq \Delta_K^{-1}, \ J(\bmu) \geq \Delta_K K \}.
$$
We say that a bandit algorithm is uniformly good if for any $\Delta_K \rightarrow 0$ at a sub-polynomial rate,
\begin{equation} \label{com}
\sup_{\bmu \in \Theta(\Delta_K)} R_N(\bmu) = o(K N^{\epsilon}) \mbox{ for all } \epsilon > 0.
\end{equation}
We show in Section 6.1 that if $\mu_k \stackrel{\rm i.i.d.}{\sim} {\rm N}(\mu_0,\sigma_0^2)$ for any $\mu_0$ real and $\sigma_0^2 > 0$,
then for $\Delta_K = (\log K)^{-\eta}$ with $\eta>\tfrac{1}{2}$,
\begin{equation} \label{Pomega}
P ( \bmu \in \Theta(\Delta_K)) \rightarrow 1.  
\end{equation}
Let $r(\bmu) = \sum_{k: \mu_k < \mu_*} \frac{2}{\mu_*-\mu_k}$.

\begin{thm} \label{thm1a}
If a bandit algorithm is uniformly good,
then for all $\Delta_K~\rightarrow~0$ at a sub-polynomial rate,
\begin{equation} \label{lower2}
\liminf_{N \rightarrow \infty} \Big[ \inf_{\bmu \in \Theta(\Delta_K)} \tfrac{R_N(\bmu)}{r(\bmu) \log N} \Big] \geq 1-\zeta. 
\end{equation}
\end{thm}

In Theorem \ref{thm1b} below for $K \rightarrow \infty$, 
for technical reasons we perturb (\ref{L1}) to 
\begin{equation} \label{Uq}
U_k^n = \bar X_{kn_k} + \sqrt{\tfrac{2 \log(n/K^{1-q})}{n_k}} \mbox{ for } q > 0.
\end{equation}
Optimality is achieved by selecting $q$ arbitrarily small,
this justifies (\ref{L1}).
For $K$ fixed,
we consider 
\begin{equation} \label{Kfixed}
U_k^n = \bar X_{kn_k} + \sqrt{\tfrac{2b_{n/K}}{n_k}},
\end{equation}
with $b_m = \log m + o(\sqrt{\log m})$ as $m \rightarrow \infty$.

We define the regret ignoring the initial allocation of one reward to each arm to be
\begin{equation} \label{Rtilde}
\wtd R_N(\bmu) = \sum_{k=1}^K (\mu_*-\mu_k) E_{\bmu}(N_k-1)^+.
\end{equation}

\begin{thm} \label{thm1b}
Consider UCB-Large as given in {\rm (\ref{Uq})} for $K \rightarrow \infty$,
or {\rm (\ref{Kfixed})} for $K$ fixed.
For $\Delta_K \rightarrow 0$ at a sub-polynomial rate, 
$$\limsup_{N \rightarrow \infty} \Big[ \sup_{\bmu \in \Theta(\Delta_K)} \tfrac{\wtd R_N(\bmu)}{r(\bmu) \log N} \Big] \leq 1-\zeta+\zeta q. 
$$
\end{thm}

Theorem \ref{thm1b} does not hold with $R_N(\bmu)$ in place of $\wtd R_N(\bmu)$.
Consider for example $\mu_1=\mu_*$ and $\mu_k = \mu_*-\log N$ for $k \geq 2$.
Here $r(\bmu) \log N = 2(K-1)$ whereas $R_N(\bmu) \geq (K-1) \log N$ due to the initial allocation of one reward to each arm under UCB-Large.

\subsection{Normal rewards with unknown and unequal variances}

Let $X_{kt} \stackrel{\rm i.i.d.}{\sim}$ N$(\mu_k,\sigma_k^2)$,
$t \geq 1$,
be the normal rewards of arm $k$.
Let $\theta_k = (\mu_k,\sigma_k^2)$ and let 
\begin{eqnarray*}
\Theta_2(\Delta_K) & = & \{ \btheta_K: \max_{k: \mu_k < \mu_*} \mu_k \leq \mu_*-\Delta_K, \max_k |\mu_k| \leq \Delta_K^{-1}, \cr
& & \quad \Delta_K \leq \sigma_k \leq \Delta_K^{-1} \mbox{ for all } k, \ J(\bmu) \geq \Delta_K K \}.
\end{eqnarray*}
A simple extension of (\ref{Pomega}) here would be to consider $\mu_k$ i.i.d. with a normal prior and $\sigma_k$ having bounded support away from 0.

Analogous to the setting of unit variance normal considered in Section~3.1,
we say that a bandit algorithm is uniformly good if for any $\Delta_K \rightarrow 0$ at a sub-polynomial rate,
\begin{equation} \label{good2}
\sup_{\btheta \in \Theta_2(\Delta_K)} R_N(\btheta) = o(KN^{\epsilon}) \mbox{ for all } \epsilon>0.
\end{equation}
Let $r(\btheta_K)$ be as given in (\ref{lowert}).

\begin{thm} \label{thm3}
If a bandit algorithm is uniformly good, 
then for all $\Delta_K \rightarrow 0$ at a sub-polynomial rate,
$$\liminf_{N \rightarrow \infty} \Big[ \inf_{\btheta \in \Theta_2(\Delta_K)} \tfrac{R_N(\btheta)}{r(\btheta) \log N} \Big] \geq 1-\zeta. 
$$
\end{thm}

As in Section 3.1, 
for technical reasons we perturb (\ref{UCBt}) for the case $K \rightarrow \infty$, 
to
\begin{equation} \label{Ut}
U_k^n = \bar X_{kn_k} + \wht \sigma_{kn_k} \sqrt{\exp(\tfrac{2 \log(n/K^{1-q})}{n_k-1})-1} \mbox{ for } q>0.
\end{equation}
Theorem \ref{thm4} below says that optimality is achieved by selecting $q$ arbitrarily small,
this justifies (\ref{UCBt}) with $b_m=\log m$.
For $K$ fixed as $N \rightarrow \infty$,
consider
\begin{eqnarray} \label{tfixed}
U_k^n & = & \bar X_{kn_k} + \wht \sigma_{kn_k} \sqrt{\exp(\tfrac{2b_{n/K}}{n_k-1})-1}, \\ \nonumber
b_m & = & \log m+ \alpha \log (1+\log m) \mbox{ for } \alpha>1.
\end{eqnarray}

\begin{thm} \label{thm4}
Consider UCB-Large as given in {\rm (\ref{Ut})} for $K \rightarrow \infty$,
or {\rm (\ref{tfixed})} for $K$ fixed.
For $\Delta_K \rightarrow 0$ at a sub-polynomial rate,
$$\limsup_{N \rightarrow \infty} \Big[ \sup_{\btheta \in \Theta_2(\Delta_K)} \tfrac{\wtd R_N(\btheta)}{r(\btheta) \log N} \Big] \leq 1-\zeta+\zeta q, 
$$
with $\wtd R_N(\btheta)$ as defined in {\rm (\ref{Rtilde})},
with $\btheta$ replacing $\bmu$.
\end{thm}

\section{UCB adjustments for polynomially decreasing gaps}

The asymptotics in Section 3 are for gaps decreasing at a sub-polynomial rate.
We extend the asymptotics here to gaps that are polynomially small.
To avoid excessive technicalities, 
we restrict to normal rewards with known variances.

Let $X_{kt} \stackrel{\rm i.i.d.}{\sim}$ N($\mu_k,1$), 
$t \geq 1$,
be the normal rewards of arm $k$.
Let $K=N^{\zeta+o(1)}$ for some $0 \leq \zeta < 1$ as $N \rightarrow \infty$,
and let $\Delta_N (=\Delta_{\alpha N})= \alpha N^{-\eta}$ for some $\alpha > 0$ and $\eta>0$.   
We consider here $\Delta_N$ instead of $\Delta_K$ (in Section 3),
so that $\Delta_N \rightarrow 0$ with $K$ fixed,
as $N \rightarrow \infty$.

In the case of polynomially decreasing gaps,
the regret bound (\ref{com}) is not achievable.
Consider for example $\mu_1 = \mu_*$ and $\mu_k = \mu_* - \Delta_N$ for $k \geq 2$.
Here $r(\bmu) = 2(K-1) \Delta_N^{-1} \log N$.
Instead for a given $\eta$,
we require a uniformly good bandit algorithm to satisfy,
instead of (\ref{com}),
\begin{equation} \label{good7}
\sup_{\bmu \in \Theta(\Delta_N)} R_N(\bmu) = O(K \Delta_N^{-1} N^{\epsilon}) \mbox{ for all } \epsilon > 0 \mbox{ and } \alpha > 0.
\end{equation}

\begin{thm} \label{thm7}
Let $0 < \eta < \tfrac{1-\zeta}{2}$.
If a bandit algorithm is such that {\rm (\ref{good7})} holds,
then for all $\alpha > 0$,
\begin{equation} \label{smalllower}
\liminf_{N \rightarrow \infty} \Big[ \tfrac{\sup_{\bmu \in \Theta(\Delta_N)} R_N(\bmu)}{2(K-1) \Delta_N^{-1} \log N} \Big] \geq 1-\zeta-2 \eta. 
\end{equation}
\end{thm}

The smaller lower bound constant in (\ref{smalllower}) compared to (\ref{lower2}),
with $1-\zeta-2 \eta$ instead of $1-\zeta$,
is due to the additional $\Delta_N^{-1}$ in the regret bound (\ref{good7}). 
To take advantage of the smaller constant (though not fully),
we consider UCB-Large as given in (\ref{Ukn}), 
with
\begin{equation} \label{bchi}
b_m = \chi \log m \mbox{ for some } \chi > 1-\tfrac{\eta}{1-\zeta} (>\tfrac{1}{2}). 
\end{equation}
The best regret guarantee given in Theorem \ref{thm8} below is for $\chi$ arbitrarily close to $1-\tfrac{\eta}{1-\zeta}$.
In practice we do not know what $\eta$ is,
and in (\ref{smalllower}) and (\ref{smallupper}) the asymptotics are for the worst-case scenarios [largest $r(\bmu)$].
Nevertheless Theorems \ref{thm7} and \ref{thm8} address why, 
in the simulations in Section~5,
numerical performances for UCB-Large are better for $\chi=\frac{1}{2}$ compared to $\chi=1$.

\begin{thm} \label{thm8}
For UCB-Large with $b_m$ as given in {\rm (\ref{bchi})},
\begin{equation} \label{smallupper}
\limsup_{N \rightarrow \infty} [\tfrac{\sup_{\bmu \in \Theta(\Delta_N)} \wtd R_N(\bmu)}{2(K-1) \Delta_N^{-1} \log N}] \leq \chi(1-\zeta). 
\end{equation}
\end{thm}

\section{Numerical studies}

We perform simulations here for normal (Examples 4 and 5) and Bernoulli (Example 6) rewards,
confirming that UCB-Large, 
which corrects for large arm sizes, 
improves upon classical UCB algorithms which don't. 
In particular UCB-Large as given in (\ref{Ukn}) with $b_m = \chi \log m$ for $\chi=0.5$ has the best performances with regrets uniformly smaller than its competitors. 

In addition to $\chi=0.5$,
we run simulations with $\chi=0.75$ and 1.
Though by (\ref{Pomega}), 
$\min_{1 \leq k \leq K} (\mu_*-\mu_k)$ is sub-polynomial when $\mu_k$ are drawn from a normal prior,
when we average the regrets over a large number of runs,
the average may be dominated by runs with polynomially small $\max_{1 \leq k \leq K} (\mu_*-\mu_k)$. 
This explains why UCB-Large with $\chi<1$,
which is better for polynomially small $\min_{1 \leq k \leq K} (\mu_*-\mu_k)$,
performs better than when $\chi=1$. 
In addition to $b_m = \chi \log m$,
we apply UCB-Large for $b_m = \log(e-1+m)-\sqrt{\log(e-1+m)}$ (labeled as $b=\log -\sqrt{\log}$), 
motivated by (\ref{Kfixed}). 
We consider $\log(e-1+m)$ instead of $\log m$ to ensure monotonicity of $b_m$ for $m \geq 1$.

In the simulations each regret is estimated by $\sum_{k=1}^K (\mu_*-\mu_k) N_k$, 
averaged over $J=10000$ simulation runs,
for $N=20000$ rewards.
Standard errors are placed after the $\pm$ sign.

\begin{table}[t]
\begin{center}
\begin{tabular}{ll|rrrr}
& & \multicolumn{4}{c}{$K$} \cr
& & 10 & 20 & 50 & 100 \cr \hline
UCB-Large & $\chi=1$ & 144$\pm$1 & 234$\pm$1 & 441$\pm$2 & 720$\pm$2 \cr
& $\chi=0.75$ & 119$\pm$2 & 193$\pm$2 & 375$\pm$3 & 624$\pm$3 \cr
& $\chi=0.5$ & 113$\pm$4 & 179$\pm$4 & 357$\pm$6 & 587$\pm$6 \cr
& $b=\log - \sqrt{\log}$ & 118$\pm$3 & 191$\pm$3 & 375$\pm$5 & 624$\pm$6 \cr \hline
UCB-Agrawal & & 176$\pm$1 & 312$\pm$1 & 650$\pm$2 & 1150$\pm$2 \cr
Thompson & & 123$\pm$1 & 213$\pm$2 & 419$\pm$2 & 706$\pm$3 \cr
\end{tabular}
\caption{ The regrets of UCB algorithms and Thompson sampling for $K$ arms. 
The rewards are normal distributed with unit variances.
The arm means are generated from a N{\rm (0,1)} prior,
and a fresh set of means is generated in each run.
We apply Thompson sampling using the correct N$(0,1)$ prior.}
\end{center}
\end{table}

\smallskip
{\sc Example} 4. 
Consider $X_{k1}, X_{k2}, \ldots$ i.i.d. N($\mu_k,1$).
We consider UCB-Agrawal with $b_n=\log n$ [see (\ref{Agr})],
as well as UCB-Large.
We also consider Thompson sampling,
assuming a N(0,1) prior for each $\mu_k$.
That is for $n \geq K$,
we generate
$$\theta_{kn} \sim \mbox{N}(\tfrac{S_{kn_k}}{n_k+1},\tfrac{1}{n_k+1}),
\quad 1 \leq k \leq K,
$$
and sample the $(n+1)$th reward from the arm $k$ maximizing $\theta_{kn}$.
This is an advantageous set-up for Thompson sampling as its prior is used for generating $\mu_k$,
that is with $\mu_k \stackrel{\rm i.i.d.}{\sim}$ N(0,1) in each run.
We see from Table~1 that the best performer is UCB-Large with $\chi=0.5$.
All the UCB-Large algorithms perform better than UCB-Agrawal.

\begin{table}[t]
\begin{center}
\begin{tabular}{ll|rrrr}
& & \multicolumn{4}{c}{$K$} \cr
& & 10 & 20 & 50 & 100 \cr \hline
UCB-Large & $\chi=1$ & 228$\pm$2 & 410$\pm$2 & 882$\pm$2 & 1596$\pm$4 \cr
& $\chi=0.75$ & 183$\pm$2 & 336$\pm$3 & 724$\pm$6 & 1293$\pm$6 \cr
& $\chi=0.5$ & 177$\pm$5 & 306$\pm$7 & 610$\pm$8 & 1070$\pm$9 \cr
& $b=\log-\sqrt{\log}$ & 179$\pm$3 & 326$\pm$6 & 664$\pm$7 & 1128$\pm$7 \cr \hline
UCB-BK & & 273$\pm$3 & 515$\pm$2 & 1222$\pm$3 & 2398$\pm$4 \cr
Thompson & & 191$\pm$6 & 313$\pm$7 & 646$\pm$7 & 1202$\pm$7  
\end{tabular}
\caption{ The regrets of UCB algorithms and Thompson sampling for $K$ arms on normal rewards with unknown and unequal variances.
The arm means are generated from {\rm N(0,1)},
the arm variances are generated from the exponential distribution with mean $1$.
We apply Thompson sampling assuming a normal-gamma prior.
}
\end{center}
\end{table}

\smallskip
{\sc Example} 5.
Consider  $X_{k1}, X_{k2}, \ldots$ i.i.d. N($\mu_k,\sigma_k^2$).
We compare UCB-Large as given in (\ref{UCBt}),
with an initial allocation of two rewards to each arm,
against UCB-BK (Burnetas and Katehakis) and Thompson sampling. 
The simulation results in Table 2 again demonstrate significant improvements to UCB with the arm-size corrections that we introduce here. 
For each run we generate $\mu_k \stackrel{\rm i.i.d.}{\sim}$ N(0,1) and $\sigma_k^2 \stackrel{\rm i.i.d.}{\sim}$ Exp(1),
the exponential distribution with mean 1.

For Thompson sampling we assume a normal-gamma prior,
generating for $n \geq 2K$,
\begin{eqnarray*}
\sigma_{kn}^{-2} & \sim & \mbox{Gamma}(1+\tfrac{n_k}{2},1+\tfrac{n_k \hat \sigma_{kn_k}^2}{2}+\tfrac{n_k \bar X_{kn_k}^2}{1+n_k}), \quad 1 \leq k \leq K, \cr
\mu_{kn}|\sigma_{kn}^{-2} & \sim & \mbox{N}(\tfrac{n_k \bar X_{kn_k}}{1+n_k}, \tfrac{\sigma_{kn}^2}{1+n_k}),
\end{eqnarray*}
and sampling the $(n+1)$th reward from the arm $k$ maximizing $\mu_{kn}$.
The best performer is UCB-Large with $\chi=0.5$,
with $b=\log-\sqrt{\log}$ and Thompson sampling both performing relatively well. 
In Thompson sampling here we do not apply the (unknown) underlying prior.

\begin{table}[t]
\begin{center}
\begin{tabular}{ll|rrrr}
& & \multicolumn{4}{c}{$K$} \cr
& & 10 & 20 & 50 & 100 \cr \hline
UCB-Large & $\chi=1$ & 61.2$\pm$0.4 &  86.6$\pm$0.5 & 138.7$\pm$0.5 & 202.5$\pm$0.7 \cr
& $\chi=0.75$ & 48.6$\pm$0.4 & 70.9$\pm$0.5 & 112.1$\pm$0.5 & 160.7$\pm$0.6 \cr
& $\chi=0.5$ & 43.4$\pm$0.8 & 60.3$\pm$0.8 & 88.4$\pm$0.7 & 133.0$\pm$0.6 \cr
& $b=\log-\sqrt{\log}$ & 46.7$\pm$0.5 & 66.3$\pm$0.7 & 96.8$\pm$0.6 & 139.6$\pm$0.6 \cr \hline
UCB-Agrawal & & 76.3$\pm$0.4 & 116.2$\pm$0.6 & 205.3$\pm$0.8 & 323.7$\pm$1.1 \cr
Thompson & & 53.3$\pm$0.4 & 79.4$\pm$0.5 & 135.4$\pm$0.6 & 207.4$\pm$0.8 \cr
\end{tabular}
\caption{ The regrets of UCB algorithms and Thompson sampling for $K$ arms. 
Rewards from arm~$k$ are Bernoulli random variables with success probabilities $p_k$,
with each $p_k$ generated from a Uniform$(0,1)$ prior. 
A fresh set of $p_k$ is generated in each run.
We apply Thompson sampling assuming the correct Uniform$(0,1)$ prior.}
\end{center}
\end{table}

\smallskip
{\sc Example} 6. 
Consider $X_{k1}, X_{k2}, \ldots$ i.i.d. Bernoulli with success probability $p_k$. 
We compare UCB-Large against UCB-Agrawal with $b_n=\log n$,
see (\ref{ruX}) and (\ref{Uk}), 
with $I_u$ as given in (\ref{LDB}). 
For Thompson sampling we assume a uniform prior for $p_k$,
that is for $n \geq K$,
we generate 
$$p_{kn} \sim \mbox{Beta}(1+S_{kn_k},1+n_k-S_{kn_k}),
\quad 1 \leq k \leq K,
$$
and sample the $(n+1)$th reward from the arm $k$ maximizing $p_{kn}$.

For each $1 \leq j \leq J$,
we generate $p_k \stackrel{\rm i.i.d.}{\sim}$ Uniform(0,1).
The arm mean $p_k$ differs with $j$.
The simulation results in Table 3 show that UCB-Large is the best performer,
its regret when $\chi=0.5$ at $K=100$ is two thirds that of Thompson sampling,
despite Thompson sampling having the advantage of applying the underlying uniform prior of the arm means.
Its regret is less than half that of UCB-Agrawal.  

\section{Proofs of (\ref{Pomega}) and Theorems \ref{thm1a} and \ref{thm1b}}

We prove (\ref{Pomega}) in Section~6.1,
Theorem \ref{thm1a} in Section 6.2 and Theorem \ref{thm1b} in Section 6.3.
Let $\phi$ denote the density and $\Phi$ the cumulative distribution of the standard normal.
Let $a_n \sim b_n$ if $\lim_{n \rightarrow \infty} \frac{a_n}{b_n} =1$ and let $\lceil \cdot \rceil$ be the least integer function.

\subsection{Proof of {\rm (\ref{Pomega})}}

Assume without loss of generality $\mu_0=0$ and $\sigma_0=1$.
Let $\xi_K$ be such that 
\begin{equation} \label{xiK}
P(\mu_* \leq \xi_K) = (2 \log K)^{-1}.
\end{equation}
Since $1-\Phi(z) \sim \tfrac{\phi(z)}{z}$ as $z \rightarrow \infty$,
$$P(\mu_* \leq \xi_K) = \{ 1-[1+o(1)] \tfrac{\phi(\xi_K)}{\xi_K}  \}^K = \exp \{ -[1+o(1)] \tfrac{K \phi(\xi_K)}{\xi_K}  \}.
$$
Replacing the above into (\ref{xiK}) leads to
\begin{equation} \label{gxi} 
\tfrac{\phi(\xi_K)}{\xi_K} \sim \tfrac{\log \log K}{K} (\Rightarrow \xi_K \sim \sqrt{2 \log K}).
\end{equation}

It follows from (\ref{xiK}),
the monotonicity of $\phi(z)$ for $z \geq 0$ and $\xi_K > \Delta_K$ for $K$ large,
that
\begin{equation} \label{gap}
P \Big( \min_{k: \mu_k < \mu_*} (\mu_*-\mu_k) \leq \Delta_K \Big) \leq K \Delta_K \phi(\xi_K-\Delta_K)+(2 \log K)^{-1}.
\end{equation}
It follows from the last relation in (\ref{gxi}) and $\Delta_K = (\log K)^{-{\eta}}$ for $\eta > \frac{1}{2}$ that $\phi(\xi_K-\Delta_K) \sim \phi(\xi_K)$, 
and therefore by (\ref{gxi}),
the probability in (\ref{gap}) goes to 0 as $K \rightarrow \infty$. 

To complete the proof,
check that
$$P(\max_{1 \leq k \leq K} |\mu_k| \geq \Delta_K^{-1}) \leq 2K \exp(-\tfrac{1}{2 \Delta_K^2}) \rightarrow 0.
$$

\subsection{Proof of Theorem {\rm \ref{thm1a}}}

Let $\lambda = \mu_*+\delta \Delta_K$ for some $\delta >0$.
For $\mu_k \leq \mu_*-\Delta_K$,
$I_{\lambda}(\mu_k) \leq (1+\delta)^2 I_{\mu_*}(\mu_k)$.
We show below that for $0 < c < 1-\zeta$,
\begin{equation} \label{show1}
E_{\bmu} N_k \geq [c+o(1)] \tfrac{\log N}{I_{\lambda}(\mu_k)} \mbox{ whenever } \mu_k \leq \mu_*-\Delta_K,
\end{equation}
with $o(1)$ uniform over $k$ and $\bmu \in \Theta(\Delta_K)$.
Theorem \ref{thm1a} follows from (\ref{show1}) by selecting $\delta$ arbitrarily small and $c$ close to $1-\zeta$.

Let $k$ be such that $\mu_k \leq \mu_* -\Delta_K$, 
and let $\blambda$ be $\bmu$ with $\lambda$ replacing $\mu_k$.
Let $a$ be such that $c < a < 1-\zeta$, 
and let
\begin{eqnarray*}
\ell_k & = & \sum_{t=1}^{N_k} \log Y_t, 
\mbox{ where } Y_t = \tfrac{f_{\mu_k}(X_{kt})}{f_{\lambda}(X_{kt})}, \cr
A_k &  = & \{ N_k < \tfrac{c \log N}{I_{\lambda}(\mu_k)}, \ell_k \leq a \log N \}, \cr
B_k &  = & \{ N_k < \tfrac{c \log N}{I_{\lambda}(\mu_k)}, \ell_k > a \log N \},
\end{eqnarray*}
with $\ell_k=0$ when $N_k=0$.
We conclude (\ref{show1}) by showing that $P_{\bmu}(A_k) \rightarrow 0$ and $P_{\bmu}(B_k) \rightarrow 0$ uniformly over $k$ and $\bmu$.

It follows from a change of measure that
\begin{equation} \label{cm}
P_{\bmu}(A_k) = E_{\blambda}(e^{\ell_k} A_k) \leq N^a P_{\blambda}(A_k).
\end{equation}
Since $1-a > \zeta$,
by the uniformly good property (\ref{com}) for the sequence $\delta \Delta_K (\rightarrow 0$ as a sub-polynomial rate), 
\begin{equation} \label{Rbound}
R_N(\blambda) = o(N^{1-a} \Delta_K). 
\end{equation}
Since $N_k = o(N)$ uniformly on $A_k$,
it follows from (\ref{cm}), 
(\ref{Rbound}) and $R_N(\blambda) \geq (\lambda-\mu_*) E_{\blambda}[(N-N_k) {\bf 1}_{A_k}]$ that
\begin{eqnarray} \label{lowerP}
[o(1) & =] & N^{a-1} \Delta_K^{-1} (\lambda-\mu_*)  E_{\blambda}[(N-N_k) {\bf 1}_{A_k}] \\ \nonumber
& \geq & [\delta+o(1)] N^a P_{\blambda}(A_k) \geq [\delta+o(1)] P_{\bmu}(A_k),
\end{eqnarray}
and so $P_{\bmu}(A_k) \rightarrow 0$.

Let $\omega=a-c$,  
$s_k = \tfrac{c \log N}{I_{\lambda}(\mu_k)}$ and $Z_t = \log Y_t-I_{\lambda}(\mu_k) [=(\mu_k-\lambda)(X_{kt}-\mu_k)]$.
It follows from the reflection principle that 
$$P_{\bmu}(B_k) \leq P_{\bmu} \Big( N_k < s_k, \sum_{t=1}^{N_k} Z_t \geq \omega \log N \Big) \leq 2[1-\Phi(\tfrac{\omega \log N}{(\lambda-\mu_k) \sqrt{s_k}})] \rightarrow 0.
$$

\subsection{Proof of Theorem {\rm \ref{thm1b}}}

Let
\begin{equation} \label{Ukr} 
U_{ks}^n = \bar X_{ks} + \sqrt{\tfrac{2 b_{nK}}{s}}, 
\end{equation}
where $b_{nK} = \log(\frac{n}{K^{1-q}})$ if $K \rightarrow \infty$,
and $b_{nK} = \log(\tfrac{n}{K})+o(\sqrt{\log(\tfrac{n}{K})})$ (as $n \rightarrow \infty$) if $K$ is fixed.
Hence $U_k^n = U_{kn_k}^n$.
Let $v = \mu_*-\delta \Delta_K$ for some $0 < \delta < 1$.
Hence $I_v(\mu_k) \geq (1-\delta)^2 I_{\mu_*}(\mu_k)$ for $\mu_k \leq \mu_*-\Delta_K$.
Let 
$$G_{ks} = \{ U_{ks}^n \geq v \mbox{ for  some } K \leq n \leq N-1 \}.
$$
Let $s_k = \tfrac{c \log N}{I_v(\mu_k)}$ for $\mu_k \leq \mu_*-\Delta_K$,
with $c>1-\zeta+\zeta q$.
We preface the proof of Theorem~\ref{thm1b} with the following lemmas,
which hold uniformly over $\bmu \in \Theta(\Delta_K)$.

\begin{lem} \label{lem2}
There exists $\beta>0$ such that 
$$\max_{k: \mu_k \leq \mu_*-\Delta_K} \sum_{s \geq s_k} P_{\mu_k}(G_{ks}) = O(N^{-\beta}).
$$
\end{lem}

\begin{lem} \label{lem3}
Let $H_n = \{ \inf_{s \geq 1} U_{\ell s}^n \leq v \}$. 
As $N \rightarrow \infty$,
\begin{equation} \label{sum1}
\sum_{n=K}^{N-1} P_{\mu_*}(H_n) = \left\{ \begin{array}{ll} O(\Delta_K^{-2} K^{1-q} \log N) & \mbox{ if } K \rightarrow \infty, \cr
o(\log N) & \mbox{ if } K \mbox{ is fixed.} \end{array} \right.
\end{equation}
\end{lem}

{\sc Proof of Theorem} \ref{thm1b}.
When there are $n$ total rewards,
an inferior arm $k$ with $s (\geq s_k)$ rewards is sampled only if either $U_{ks}^n \geq v$, 
or $\inf_{t \geq 1} U_{\ell t}^n \leq v$ for an optimal arm $\ell$.
Hence
\begin{eqnarray} \label{Nks}
& & \sum_{k: \mu_k < \mu_*} (N_k-s_k-1)^+ \leq \sum_{k: \mu_k < \mu_*} (N_k- \lceil s_k \rceil)^+ \\ \nonumber
& \leq & \sum_{n=K}^{N-1} {\bf 1}_{\{ {\rm arm } \ k \ {\rm sampled \ with } \ U_{kn_k}^n \geq v \ 
{\rm for \ some } \ k \ {\rm such \ that } \ n_k \geq s_k \ {\rm or } \ U_{ln_l}^n \leq v \}} \\ \nonumber
& \leq & \sum_{k: \mu_k < \mu_*} \sum_{s \geq s_k} {\bf 1}_{G_{ks}} + \sum_{n=K}^{N-1} {\bf 1}_{H_n}.
\end{eqnarray}
By Lemmas \ref{lem2} and \ref{lem3},
the right-hand side of (\ref{Nks}) is $o(K \Delta_K^2 \log N)$ after taking expectations.
Since $(N_k-1)^+ \leq (N_k-s_k-1)^+ + s_k$,
Theorem \ref{thm1b} follows from selecting $\delta$ arbitrarily small and $c$ close to $1-\zeta+\zeta q$.
$\wbox$ 

\smallskip
{\sc Proof of Lemma} \ref{lem2}.
Since $c>1-\zeta+\zeta q$,
there exists $\epsilon>0$ be such that $\omega := (\tfrac{1-\zeta+\zeta q+2 \epsilon}{c})^{\frac{1}{2}} < 1$.
By (\ref{Ukr}),
if $U_{ks}^n \geq v$,
then for $N$ large,
\begin{eqnarray*} 
(Z_{ks} :=)  \sqrt{s}(\bar X_{ks}-\mu_k) & \geq & \sqrt{s}(v-\mu_k)-\sqrt{2b_{nK}} \cr
& \geq & \sqrt{s} (v-\mu_k)-\sqrt{2(1-\zeta+\zeta q+\epsilon) \log N}. 
\end{eqnarray*}
For $K$ fixed,
$\zeta=\zeta q=0$,
so the above inequalities still hold.
Hence under $G_{ks}$ for $s \geq s_k$,
\begin{eqnarray*}
Z_{ks} & \geq & \sqrt{s}(1-\omega)(v-\mu_k)+\sqrt{s} \omega(v-\mu_k) -\sqrt{2(1-\zeta+\zeta q+\epsilon) \log N} \cr
& \geq & \sqrt{s}(1-\omega)(v-\mu_k)+\sqrt{2 \beta \log N},
\end{eqnarray*}
where $\beta = (\omega \sqrt{c}-\sqrt{1-\zeta+\zeta q+\epsilon})^2$.
Hence 
$$\sum_{s \geq s_k} P_{\mu_k}(G_{ks}) \leq N^{-\beta} \sum_{s \geq s_k} e^{-\frac{s(v-\mu_k)^2(1-\omega)^2}{2}},
$$
and Lemma \ref{lem2} holds. 
$\wbox$ 

\smallskip
{\sc Proof of Lemma} \ref{lem3}.
Let $Z_s = \sqrt{s}(\mu_*-\bar X_{\ell s})$.
By (\ref{Ukr}) and $\mu_*-v=\delta \Delta_K$,
\begin{eqnarray} \label{eq2} 
P_{\mu_*}(U_{\ell s}^n \leq v) & = & P(Z_s \geq \delta \Delta_K \sqrt{s}+\sqrt{2b_{nK}}) \\ \nonumber
& \leq & \exp(-\tfrac{s \delta^2 \Delta_K^2}{2}-b_{nK}-\delta \Delta_K \sqrt{2b_{nK}}).
\end{eqnarray}
Let $C = \sup_{y > 0} \tfrac{y^2 e^{-y^2/2}}{1-e^{-y^2/2}} (< \infty)$.
By (\ref{eq2}),
\begin{eqnarray} \label{PHn}
P_{\mu_*}(H_n) & \leq & \sum_{s=1}^{\infty} P_{\mu_*}(U_{\ell s}^n \leq v) \\ \nonumber
& \leq & 2C (\delta \Delta_K)^{-2} \exp(-b_{nK}-\delta \Delta_K \sqrt{2b_{nK}}).
\end{eqnarray}

For $K \rightarrow \infty$ with $b_{nK} = \log(\tfrac{n}{K^{1-q}})$,
$$\sum_{n=K}^{N-1} e^{-b_{nK}} \leq K^{1-q} \log N,
$$
and (\ref{sum1}) follows from (\ref{PHn}).

Consider next $K$ fixed with $b_{nK} = b_{n/K}$, 
where 
\begin{equation} \label{bm6}
b_m = \log m+o(\sqrt{\log m}).
\end{equation}
By (\ref{PHn}) it suffices to show that for any $\omega > 0$,
\begin{equation} \label{omega}
\sum_{n=K}^{N-1} \exp(-b_{nK}-\delta \Delta_K \sqrt{2b_{nK}}) \leq \omega \log N \mbox{ for } N \mbox{ large.}
\end{equation}
Let $\tau>0$ be such that $e^{-\tau} K < \omega$,
and note that by (\ref{bm6}),
there exists positive integer $m_{\tau}$ such that
$$b_m+\delta \Delta_K \sqrt{2b_m} \geq \log m+\tau \mbox{ for } m \geq m_{\tau}.
$$
Hence by the monotonicity of $b_m$, 
\begin{eqnarray*} 
\sum_{n=K}^{N-1} \exp(-b_{nK}-\delta \Delta_K \sqrt{2b_{nK}}) & \leq & Km_{\tau} e^{-b_1}+e^{-\tau} \sum_{n=Km_{\tau}+1}^{N-1} \tfrac{K}{n} \cr
& \leq & Km_{\tau} e^{-b_1}+e^{-\tau} K \log N,
\end{eqnarray*}
and (\ref{omega}) follows from $e^{-\tau} K < \omega$. 
$\wbox$

\section{Proofs of Theorems \ref{thm3} and \ref{thm4}}

We prove Theorems \ref{thm3} and \ref{thm4} in Sections~7.1 and 7.2 respectively.

\subsection{Proof of Theorem \ref{thm3}}

For a given $\btheta \in \Theta(\Delta_K)$,
let $\lambda = (\mu_{\lambda},\sigma^2_{\lambda})$,
with $\mu_{\lambda}=\mu_*+\delta \Delta_K$ for some $\delta >0$ and $\sigma^2_{\lambda} = \sigma_k^2 + (\mu_{\lambda}-\mu_k)^2$.
For $z > 0$, 
\begin{equation} \label{kappaz}
\log (1+\kappa z) \left\{ \begin{array}{ll} \leq \kappa \log(1+z) & \mbox{ if } \kappa >1, \cr
\geq \kappa \log(1+z) & \mbox{ if } \kappa <1. \end{array} \right.
\end{equation} 
To show (\ref{kappaz}),
check that equality holds at $z=0$, 
and that the first derivatives with respect to $z$ follow the inequalities. 

It follows from (\ref{kappaz}) that if $\mu_k \leq \mu_*-\Delta_K$,
then $M(\tfrac{\mu_{\lambda}-\mu_k}{\sigma_k}) \leq (1+\delta)^2 M(\tfrac{\mu_*-\mu_k}{\sigma_k})$.
We show below that for $0 < c <1-\zeta$,
\begin{equation} \label{EtN} 
E_{\btheta} N_k \geq [c+o(1)] \tfrac{\log N}{M(\frac{\mu_{\lambda}-\mu_k}{\sigma_k})} \mbox{ whenever } \mu_k \leq \mu_*-\Delta_K, 
\end{equation}
with $o(1)$ uniform over $k$ and $\btheta \in \Theta(\Delta_K)$.
Theorem \ref{thm3} follows from (\ref{EtN}) by selecting $\delta$ arbitrarily small and $c$ close to $1-\zeta$.

Let $k$ be such that $\mu_k \leq \mu_*-\Delta_K$,
and let $\blambda$ be $\btheta$ with $\lambda$ replacing $\theta_k$.
Let $a$ be such that $c < a < 1-\zeta$, 
and let
\begin{eqnarray*}
\ell_k & = & \sum_{t=1}^{N_k} \log Y_t, \mbox{ where } Y_t = \tfrac{f_{\theta_k}(X_{kt})}{f_{\lambda}(X_{kt})}, \cr
A_k & = & \{ N_k < \tfrac{c \log N}{M(\frac{\mu_{\lambda}-\mu_k}{\sigma_k})}, \ \ell_k \leq a \log N \}, \cr
B_k & = & \{ N_k < \tfrac{c \log N}{M(\frac{\mu_{\lambda}-\mu_k}{\sigma_k})}, \ \ell_k > a \log N \}.
\end{eqnarray*}
It follows from $P_{\btheta}(A_k) = E_{\blambda}(e^{\ell_k} A_k) \leq N^a P_{\blambda}(A_k)$, 
the uniformly good property (\ref{good2}) and the computations in (\ref{lowerP}) that $P_{\btheta}(A_k) \rightarrow 0$,
uniformly over $k$ and $\btheta$,
and so (\ref{EtN}) follows from $P_{\btheta}(B_k) \rightarrow 0$.

Let $\omega=a-c$,
$s_k = \tfrac{c \log N}{M(\frac{\mu_{\lambda}-\mu_k}{\sigma_k})}$,
$V_t = \log Y_t - M(\tfrac{\mu_{\lambda}-\mu_k}{\sigma_k})$ and check that 
\begin{equation} \label{PtB}
P_{\btheta}(B_k) \leq P_{\btheta} \Big( N_k < s_k, 
\sum_{t=1}^{N_k} V_t \geq \omega \log N \Big).
\end{equation}
Let $Z_t = \frac{X_{kt}-\mu_k}{\sigma_k}$ and $d_k = \frac{\mu_{\lambda}-\mu_k}{\sigma_k}$.
Hence $\frac{X_{kt}-\mu_{\lambda}}{\sigma_{\lambda}} = \frac{\sigma_k}{\sigma_{\lambda}}(Z_t-d_k)$, 
$\frac{\sigma_k^2}{\sigma_{\lambda}^2} = \frac{1}{1+d_k^2}$ and 
\begin{eqnarray} \label{Vt}
V_t & = & \tfrac{(X_{kt}-\mu_{\lambda})^2}{2 \sigma_{\lambda}^2}-\tfrac{(X_{kt}-\mu_k)^2}{2 \sigma_k^2} +\tfrac{1}{2} \log (\tfrac{\sigma_{\lambda}^2}{\sigma_k^2}) 
-\tfrac{1}{2} \log[1+(\tfrac{\mu_{\lambda}-\mu_k}{\sigma_k})^2] \\ \nonumber
& = & \tfrac{1}{2}(\tfrac{\sigma_k^2}{\sigma_{\lambda}^2}-1) Z_t^2 - \tfrac{d_k \sigma_k^2}{\sigma_{\lambda}^2} Z_t + \tfrac{\sigma_k^2 d_k^2}{2 \sigma_{\lambda}^2} \\ \nonumber
& = & \tfrac{d_k^2}{2(1+d_k^2)}(1-Z_t^2)-\tfrac{d_k}{1+d_k^2} Z_t.
\end{eqnarray}
Let $C_1 = \sup_{d>0} \tfrac{d^2/(1+d^2)}{M(d)} (<\infty)$ and $C_2 = \sup_{d>0} \tfrac{d/(1+d^2)}{\sqrt{M(d)}} (<\infty)$.
By (\ref{Vt}) and,
for $0 < x < 1$,
$$P(\chi_n^2/n \leq 1-x) \leq \exp \{ n[x+\log(1-x)] \} \leq \exp(-\tfrac{nx^2}{2}), 
$$
where $\chi^2_n$ is a $\chi^2$-random variable with $n$ degrees of freedom,
\begin{eqnarray} \label{Pv}
& & P_{\theta_k} \Big( \sum_{t=1}^n V_t \geq \omega \log N \Big) \\ \nonumber
& \leq & P_{\theta_k} \Big( \sum_{t=1}^n Z_t^2 \leq n-\tfrac{\omega \log N}{C_1 M(d_k)} \Big) + P_{\theta_k} \Big( \sum_{t=1}^n Z_t \leq -\tfrac{\omega \log N}{2C_2 \sqrt{M(d_k)}} \Big) \\ \nonumber
& \leq & \exp[-\tfrac{(\omega \log N)^2}{2n C_1^2 M^2(d_k)}]+\exp[-\tfrac{(\omega \log N)^2}{8n C_2^2 M(d_k)}].
\end{eqnarray}
Since $s_k = \frac{c \log N}{M(d_k)}$ and $M(d_k) = o(\log N)$,
summing (\ref{Pv}) over $1 \leq n < s_k$ and substituting into (\ref{PtB}) leads to $P_{\btheta}(B_k) \rightarrow 0$.

\subsection{Proof of Theorem \ref{thm4}}

Let
\begin{eqnarray} \label{U7}
& & U_{ks}^n = \bar X_{ks} + \wht \sigma_{ks} \sqrt{\exp(\tfrac{2b_{nK}}{s-1})-1} \\ \nonumber
& & [\Rightarrow (s-1) M(\tfrac{U_{ks}^n-\bar X_{ks}}{\hat \sigma_{ks}}) = b_{nK}].
\end{eqnarray}
Let $v = \mu_*-\delta \Delta_K$ for some $0 < \delta < 1$. 
By (\ref{kappaz}),
$$M(\tfrac{v-\mu_k}{\sigma_k}) \geq (1-\delta)^2 M(\tfrac{\mu_*-\mu_k}{\sigma_k}) \mbox{ for } \mu_k \leq \mu_*-\Delta_K.
$$
Let
$$G_{ks} = \{ U_{ks}^n \geq v \mbox{ for some } K \leq n \leq N-1 \}.
$$
Let $s_k = \tfrac{c \log N}{M(\frac{v-\mu_k}{\sigma_k})}+1$ for $\mu_k \leq \mu_*-\Delta_K$, 
with $c>1-\zeta+\zeta q$.
Theorem \ref{thm4} follows from Lemmas \ref{lem3a} and \ref{lem4a} below [which hold uniformly over $\btheta \in \Theta_2(\Delta_K)$] and (\ref{Nks}),
with $\delta$ selected arbitrarily small and $c$ close to $1-\zeta+\zeta q$.
 
\begin{lem} \label{lem3a}
As $N \rightarrow \infty$, 
$$\max_{k: \mu_k \leq \mu_*-\Delta_K} \sum_{s \geq s_k} P_{\theta_k}(G_{ks}) \rightarrow 0. 
$$
\end{lem}

\begin{lem} \label{lem4a}
Let $H_n = \{ \inf_{s \geq 2} U_{\ell s}^n \leq v \}$.
As $N \rightarrow \infty$,
$$\sum_{n=K}^{N-1} P_{\theta_*}(H_n) = \left\{ \begin{array}{ll} O( \Delta_K^{-6} K^{1-q} (\log N)^2 ) & \mbox{ if } K \rightarrow \infty, \cr
o(\log N) & \mbox{ if } K \mbox{ is fixed.} \end{array} \right.
$$
\end{lem} 

{\sc Proof of Lemma} \ref{lem3a}.
Since $c > 1-\zeta+\zeta q$,
there exists $\epsilon>0$ such that $\omega := (\frac{1-\zeta+\zeta q+\epsilon}{c})^{\frac{1}{4}} < 1$.
Since $b_{nK} \leq b_{NK} \leq (1-\zeta+\zeta q+\epsilon) \log N$ for $N$ large,
by (\ref{U7}),
under $G_{ks}$ for $s \geq s_k$,
either $v \leq \bar X_{ks}$ or 
\begin{equation} \label{Mv}
M(\tfrac{v-\bar X_{ks}}{\hat \sigma_{ks}}) \leq (\tfrac{1-\zeta+\zeta q+\epsilon}{s_k-1}) \log N = \omega^4 M(\tfrac{v-\mu_k}{\sigma_k}).
\end{equation}
By (\ref{kappaz}) and (\ref{Mv}),
$(\tfrac{v-\bar X_{ks}}{\hat \sigma_{ks}})^2 \leq \omega^4 (\tfrac{v-\mu_k}{\sigma_k})^2$.
Since $v > \mu_k$, 
this implies that under $G_{ks}$ for $s \geq s_k$,
either 
\begin{equation} \label{either}
v-\bar X_{ks} \leq \omega (v-\mu_k) \mbox{ or } \wht \sigma_{ks} \geq \omega^{-1} \sigma_k.
\end{equation}
Let $\eta = \omega^{-2}-1-\log \omega^{-2} (>0)$.
We conclude from (\ref{either}) and
$$P(\chi_{s-1}^2/(s-1) \geq x) \leq \exp[-(s-1)(x-1-\log x)] \mbox{ for } x > 1,
$$
that for $s \geq s_k$,
$$P_{\theta_k}(G_{ks}) \leq \exp[-\tfrac{s(1-\omega)^2(v-\mu_k)^2}{2 \sigma_k^2}] + \exp[-(s-1) \eta],
$$
and Lemma \ref{lem3a} holds because $M(\frac{v-\mu_k}{\sigma_k}) = o(\log N)$ under $\Theta_2(\Delta_K)$ and $\inf_{d > 0} \frac{d^2}{M(d)} > 0$.
$\wbox$

\smallskip
{\sc Proof of Lemma} \ref{lem4a}.
Let
\begin{equation} \label{HJ}
H_{ns} = \{ U_{\ell s}^n \leq v, \wht \sigma_{\ell s}^2 \geq \sigma_{\ell}^2 \}, \quad J_{ns} = \{ U_{\ell s}^n \leq v, \wht \sigma_{\ell s}^2 \leq \sigma_{\ell}^2 \}.
\end{equation}
Let $Z_{\ell s} = \frac{\sqrt{s}(\mu_*-\bar X_{\ell s})}{\sigma_{\ell}}$.
If $U_{\ell s}^n \leq v$,
then by (\ref{U7}),
$$(s-1) M(\tfrac{v-\bar X_{\ell s}}{\hat \sigma_{\ell s}}) \geq b_{nK}, \quad v \geq \bar X_{\ell s},
$$
with $b_{nK} = \log(\tfrac{n}{K^{1-q}})$ if $K \rightarrow \infty$, 
and $b_{nK} = \log(\tfrac{n}{K})+\alpha \log [1+\log(\tfrac{n}{K})]$ for $\alpha >1$ if $K$ is fixed.
Hence under $H_{ns}$,
by (\ref{U7}),
in view that $\sigma_{\ell}^{-1} \geq \Delta_K$ and $e^y-1 \geq y$,
\begin{eqnarray*} 
Z_{\ell s} = \tfrac{\sqrt{s}(v-\bar X_{\ell s}+\delta \Delta_K)}{\sigma_{\ell}} & \geq & s^{\frac{1}{2}} \{[\exp(\tfrac{2b_{nK}}{s-1})-1]^{\frac{1}{2}}+\tfrac{\delta \Delta_K}{\sigma_{\ell}} \} \cr
& \geq & (2b_{nK})^{\frac{1}{2}}+s^{\frac{1}{2}} \delta \Delta_K^2,
\end{eqnarray*}
therefore 
\begin{eqnarray} \label{sumH}
\sum_{n=K}^{N-1} \sum_{s=2}^{\infty} P_{\theta_*}(H_{ns}) 
& \leq & \Big( \sum_{n=K}^{N-1} e^{-b_{nK}} \Big) \Big( \sum_{s=2}^{\infty} e^{-\frac{s(\delta \Delta_K^2)^2}{2}} \Big) \\ \nonumber
& = & \left\{ \begin{array}{ll} O(\Delta_K^{-4} K^{1-q} \log N) & \mbox{ if } K \rightarrow \infty, \cr
o(\log N) & \mbox{ if } K \mbox{ is fixed.} \end{array} \right.
\end{eqnarray}

It remains to show analogous bounds with $J_{ns}$ in place of $H_{ns}$.
Under $J_{ns}$,
$\wht \sigma_{\ell s}^2 \leq \sigma_{\ell}^2$ and $\bar X_{\ell s} < U_{\ell s}^n \leq v < \mu_*$,
hence by (\ref{U7}) and $\frac{\mu_*-v}{\sigma_{\ell}} \geq \delta \Delta_K^2$,
\begin{equation} \label{3b2}
\tfrac{(\mu_*-\bar X_{\ell s})^2}{\hat \sigma_{\ell s}^2} \geq \tfrac{(\mu_*-v)^2}{\hat \sigma_{\ell s}^2}+\tfrac{(U_{\ell s}^n-\bar X_{\ell s})^2}{\hat \sigma_{\ell s}^2} 
\geq \delta^2 \Delta_K^4 + \exp(\tfrac{2b_{nK}}{s-1})-1 (:= \kappa_{ns}).
\end{equation}
It follows from (\ref{U7}),
(\ref{HJ}) and (\ref{3b2}) that
\begin{eqnarray} \label{3b3}
J_{ns} & \subset & \{ \tfrac{(\mu_*-\bar X_{\ell s})^2}{\hat \sigma_{\ell s}^2} \geq \kappa_{ns}, \bar X_{\ell s} < \mu_* \} \\ \nonumber
& = & \{ T_{s-1} \geq t_{ns} \},  
\end{eqnarray}
where $t_{ns} = [(s-1) \kappa_{ns}]^{\frac{1}{2}}$ and $T_{s-1} = \tfrac{\sqrt{s-1}(\mu_*-\bar X_{\ell s})}{\hat \sigma_{\ell s}}$.

Under $P_{\theta_*}$,
$T_{s-1}$ has a $t$-distribution with $(s-1)$ degrees of freedom.
Hence Lemma \ref{lem4a} follows from (\ref{sumH}) and 
\begin{eqnarray} \label{App}
& & \sum_{n=K}^{N-1} \sum_{s=2}^{\infty} P(T_{s-1} \geq t_{ns}) \\ \nonumber
& = & \left\{ \begin{array}{ll} O(\Delta_K^{-6} K^{1-q} (\log N)^2) & \mbox{ if } K \rightarrow \infty, \cr
o(\log N) & \mbox{ if } K \mbox{ is fixed.} \end{array} \right.
\end{eqnarray}
We show (\ref{App}) in Appendix A.
$\wbox$

\section{Proofs of Theorems \ref{thm7} and \ref{thm8}}

We prove Theorems \ref{thm7} and \ref{thm8} in Sections 8.1 and 8.2 respectively.

\subsection{Proof of Theorem \ref{thm7}}

Consider $\bmu$ such that $\mu_1=\mu_*$ and $\mu_k = \mu_*-\Delta_N$ for $k \geq	2$.
Let $\lambda = \mu_*+\delta \Delta_N$ for some $\delta>0$ and note that $I_{\lambda}(\mu_k) = (1+\delta)^2 I_{\mu_*}(\mu_k)$.
We show below that for $0 < c < 1-\zeta-2 \eta$,
\begin{equation} \label{N7}
E_{\bmu} N_k \geq [c+o(1)] \tfrac{\log N}{I_{\lambda}(\mu_k)}. 
\end{equation}
Theorem \ref{thm7} follows from (\ref{N7}) by selecting $\delta$ arbitrarily small and $c$ close to $1-\zeta-2 \eta$.

Let $\blambda$ be $\bmu$ with $\lambda$ replacing $\mu_k$.
Let $a$ and $\epsilon$ be such that $c < a < a+\epsilon < 1-\zeta -2 \eta$, 
and let
\begin{eqnarray*}
\ell_k & = & \sum_{t=1}^{N_k} \log Y_t, \mbox{ where } Y_t = \tfrac{f_{\mu_k}(X_{kt})}{f_{\lambda}(X_{kt})}, \cr
A_k & = & \{ N_k < \tfrac{c \log N}{I_{\lambda}(\mu_k)}, \ell_k \leq a \log N \}, \cr
B_k & = & \{ N_k < \tfrac{c \log N}{I_{\lambda}(\mu_k)}, \ell_k > a \log N \},
\end{eqnarray*}
with $\ell_k=0$ when $N_k=0$.
The inequality (\ref{N7}) follows from $P_{\bmu}(A_k) \rightarrow 0$ and $P_{\bmu}(B_k) \rightarrow 0$.

It follows from the uniformly good property (\ref{good7}) for the sequence $\delta \Delta_N$, 
$P_{\bmu}(A_k) = E_{\blambda}(e^{\ell_k} A_k) \leq N^a P_{\blambda}(A_k)$ and $N_k = o(N)$ uniformly on $A_k$ that
\begin{eqnarray*}
[O(K \Delta_N^{-1} N^{\epsilon}) =] R_N(\blambda) & \geq & (\lambda-\mu_*) E_{\blambda}[(N-N_k) {\bf 1}_{A_k}] \cr
& \geq & [\delta+o(1)] \Delta_N N^{1-a} P_{\bmu}(A_k),
\end{eqnarray*}
hence $P_{\bmu}(A_k) = O(K \Delta_N^{-2} N^{a+\epsilon-1}) \rightarrow 0$.

Let $\omega=a-c$,
$s_k = \tfrac{c \log N}{I_{\lambda}(\mu_k)}$ and 
$$Z_t = \log Y_t-I_{\lambda}(\mu_k) = (\mu_k-\lambda)(X_{kt}-\mu_k).
$$
It follows from the reflection principle that 
\begin{eqnarray*}
P_{\bmu}(B_k) & \leq & P_{\bmu} \Big( \sum_{t=1}^s Z_t \geq \omega \log N \mbox{ for some } s < s_k \Big) \cr
& \leq & 2[1-\Phi(\tfrac{\omega \log N}{(\lambda-\mu_k) \sqrt{s_k}})] \leq 2 \exp(-\tfrac{\omega^2 \log N}{4c})] \rightarrow 0. 
\end{eqnarray*}

\subsection{Proof of Theorem \ref{thm8}}

Let $0 < \delta < 1$ to be further specified. 
For $k$ such that $\mu_k \leq \mu_*-\Delta_N$,
let $v_k = \mu_*-\delta (\mu_*-\mu_k)$. 
Hence $I_{v_k}(\mu_k) = (1-\delta)^2 I_{\mu_*}(\mu_k)$.
Let 
$$G_{ks} = \{ U_{ks}^n \geq v_k \mbox{ for some } K \leq n \leq N-1 \},
$$
and let $s_k = \tfrac{c \log N}{I_{v_k}(\mu_k)}$ for $c>\chi(1-\zeta)$.
We preface the proof of Theorem \ref{thm8} with the following lemmas. 

\begin{lem} \label{lem7}
There exists $\beta>0$ and $C > 0$ such that for $N$ large,
$$\sum_{s \geq s_k} P_{\mu_k}(G_{ks}) \leq C N^{-\beta}(\mu_*-\mu_k)^{-2}.
$$
\end{lem}

\begin{lem} \label{lem8}
Let $\gamma<\chi$.
There exists $C_{\gamma} > 0$ such that
$$P_{\mu_*}(\inf_{t \geq 1} U_{\ell t}^n \leq \mu_*-d) \leq C_{\gamma} d^{-1} (\tfrac{K}{n})^{\gamma}, 
$$
for $d>0$ and $K \leq n \leq N-1$.
\end{lem}

{\sc Proof of Theorem} \ref{thm8}.
We show below that
\begin{equation} \label{show8a}
\sup_{\bmu \in \Theta(\Delta_N)} \wtd R_N(\bmu) \leq [2 \chi(1-\zeta)+o(1)] (K-1) \Delta_N^{-1} \log N.
\end{equation}

When there are $n$ total rewards,
an inferior arm $k$ with $s$ rewards is sampled only if $U_{ks}^n \geq v_k$ or $\inf_{t \geq 1} U_{\ell t}^n \leq v_k$ for an optimal arm $\ell$.
Let
$$D_h = \{ k: e^h \Delta_N \leq \mu_* - \mu_k \leq e^{h+1} \Delta_N \}.
$$
Analogous to (\ref{Nks}),
since $v_k \leq \mu_* - e^h \delta \Delta_N$ for $k \in D_h$,
\begin{eqnarray} \label{expand}
& & \sum_{k \in D_h} (\mu_*-\mu_k) (N_k-s_k-1)^+ \leq \sum_{k \in D_h} \sum_{s \geq s_k} (\mu_*-\mu_k) {\bf 1}_{G_{ks}} \\ \nonumber
& & \quad +\sum_{n=K}^{N-1} e^{h+1} \Delta_N {\bf 1}_{\{ \inf_{t \geq 1} U_{\ell t}^n \leq \mu_*-e^h \delta \Delta_N \}}.
\end{eqnarray}
The complication in (\ref{expand}) compared to (\ref{Nks}) is needed due to the wider range of $\mu_* - \mu_k$ that we consider here.

Since $\eta > (1-\chi)(1-\zeta)$,
we can find $\gamma < \chi$ such that $\eta > (1-\gamma)(1-\zeta)$.
Since $(N_k-1)^+ \leq (N_k-s_k-1)^+ + s_k$,
by Lemmas \ref{lem7} and \ref{lem8}, 
summing (\ref{expand}) over $0 \leq h \leq h_N[:=\log (2 \Delta_N^{-2})]$ and taking expectations,
\begin{eqnarray*}
& & \wtd R_N(\bmu) - \sum_{k: \mu_k < \mu_*} (\mu_*-\mu_k) s_k \cr
& \leq & \sum_{k: \mu_k < \mu_*} (\mu_*-\mu_k) \sum_{s \geq s_k} P_{\mu_k}(G_{ks}) \cr
& & \quad + \sum_{0 \leq h \leq h_N} \Big[ e^{h+1} \Delta_N \sum_{n=K}^{N-1} P_{\mu_*}(\inf_{t \geq 1} U_{\ell t}^n \leq \mu_*-e^h \delta \Delta_N) \Big] \cr
& \leq & C KN^{-\beta} \Delta_N^{-1} + \delta^{-1} e C_{\gamma} (h_N+1) \sum_{n=K}^{N-1}(\tfrac{K}{n})^{\gamma} \cr
& = & O(KN^{-\beta} \Delta_N^{-1}+KN^{(1-\gamma)(1-\zeta)+o(1)}) = o(K \Delta_N^{-1}),
\end{eqnarray*}
in view that $\Delta_N^{-1} = \alpha^{-1} N^{\eta}$,
and (\ref{show8a}) follows from selecting $c$ arbitrarily close to $\chi(1-\zeta)$ and $\delta$ close to 0.
$\wbox$

\smallskip
{\sc Proof of Lemma} \ref{lem7}.
Since $c>\chi(1-\zeta)$,
there exists $\epsilon>0$ such that $\omega := (\tfrac{\chi (1-\zeta)+2 \epsilon}{c})^{\frac{1}{2}} < 1$.
If for some $s \geq s_k$ and $K \leq n \leq N-1$,
$$U_{ks}^n \Big(= \bar X_{ks} + \sqrt{\tfrac{2 \chi \log(n/K)}{s}} \Big) \geq v_k, 
$$
then for $N$ large,
\begin{eqnarray*}
\sqrt{s}(\bar X_{ks}-\mu_k) & \geq & \sqrt{s}(v_k-\mu_k) -\sqrt{2[\chi(1-\zeta)+\epsilon] \log N} \cr
& \geq & \sqrt{s}(1-\omega)(v_k-\mu_k)+\sqrt{2 \beta \log N},
\end{eqnarray*}
where $\beta = (\omega \sqrt{c}-\sqrt{\chi(1-\zeta)+\epsilon})^2$.
Hence
$$P_{\mu_k}(G_{ks}) \leq N^{-\beta} e^{-\frac{s(1-\omega)^2(v_k-\mu_k)^2}{2}},
$$
and Lemma \ref{lem7} holds with $C = \sup_{x > 0} (\tfrac{x^2 e^{-(1-\omega)^2 x^2/2}}{1-e^{-(1-\omega)^2 x^2/2}})$.
$\wbox$

\smallskip
{\sc Proof of Lemma} \ref{lem8}.
Let integer $j_0$ be such that $(\tfrac{j_0}{j_0+1}) \chi \geq \gamma$,
and let $H_j^n = \{ \inf_{jd^{-1} < t \leq (j+1)d^{-1}} U_{\ell t}^n \leq \mu_*-d \}$.
Let $W_t = t(\mu_*-\bar X_{\ell t})$. 
It follows from the reflection principle that for $j \geq j_0$,
\begin{eqnarray} \label{PHj}
& & P_{\mu_*}(H_j^n) \\ \nonumber
& = & P_{\mu_*}(W_t \geq dt+\sqrt{2t \chi \log(n/K)} \mbox{ for some } jd^{-1} < t \leq (j+1)d^{-1}) \\ \nonumber
& \leq & P_{\mu_*}(\max_{t \leq (j+1)d^{-1}} W_t \geq j + \sqrt{2jd^{-1} \chi \log(n/K)}) \\ \nonumber
& \leq & 2(\tfrac{K}{n})^{\gamma} \exp(-\tfrac{dj^2}{2(j+1)}) \leq 2(\tfrac{K}{n})^{\gamma} \exp(-\tfrac{dj}{4}).
\end{eqnarray}
For $t \leq j_0 d^{-1}$,
$$P_{\mu_*}(U_{\ell t}^n \leq \mu_*-d) \leq P_{\mu_*}(U_{\ell t}^n \leq \mu_*) \leq (\tfrac{K}{n})^{\chi} \leq (\tfrac{K}{n})^{\gamma},
$$
and therefore by (\ref{PHj}),
Lemma \ref{lem8} holds for $C_{\gamma} = j_0 + 8(\sup_{x>0} \tfrac{xe^{-x}}{1-e^{-x}})$.
$\wbox$

\begin{appendix}
\section{Proof of (\ref{App})}

The $t$-distribution with $(s-1)$ degrees of freedom has density
$$g_{s-1}(t) = C_s (1+\tfrac{t^2}{s-1})^{-\frac{s}{2}}, \mbox{ where } C_s = \tfrac{\Gamma(\frac{s}{2})}{\sqrt{(s-1) \pi} \Gamma(\frac{s-1}{2})}. 
$$
Let $\xi = \delta \Delta_K^2$.
In view that $\frac{z+\xi^2}{z+\xi^2-1} \leq 1+\xi^{-2}$ for $z =\exp(\frac{2b_{nK}}{s-1}) (\geq 1)$
[so $\kappa_{ns}^{-1} \leq (1+\xi^{-2})(1+\kappa_{ns})^{-1}$ and hence $t_{ns}^{-2} \leq (1+\xi^{-2})(s-1)^{-1}(1+\kappa_{ns})^{-1}$],
\begin{eqnarray} \label{3b5}
P(T_{s-1} \geq t_{ns}) & \leq & C_s t_{ns}^{-1} \int_{t_{ns}}^{\infty} \tfrac{t}{(1+\frac{t^2}{s-1})^{\frac{s}{2}}} dt \\ \nonumber
& = & \left\{ \begin{array}{ll} \frac{1}{2} C_2 t_{n2}^{-1} \log(1+\kappa_{n2}) & \mbox{ if } s=2, \cr
\tfrac{s-1}{s-2} C_s t_{ns}^{-1} (1+\kappa_{ns})^{-\frac{s}{2}+1} & \mbox{ if } s>2, \end{array} \right. \\ \nonumber
& \leq & C(1+\xi^{-2})^{\frac{1}{2}} y_{ns},
\end{eqnarray}
where $C = \sup_{s \geq 2} 2C_s (<\infty$ because $C_s \rightarrow \frac{1}{\sqrt{2 \pi}}$ as $s \rightarrow \infty$) and 
\begin{eqnarray} \label{yn2}
y_{n2} & = & \log[\exp(2b_{nK})+\xi^2]/[\exp(2b_{nK})+\xi^2]^{\frac{1}{2}}, \\ \label{yns}
y_{ns} & = & (s-1)^{-\frac{1}{2}} [\exp(\tfrac{2b_{nK}}{s-1})+\xi^2]^{-\frac{(s-1)}{2}}, \quad s \geq 3.
\end{eqnarray}
By (\ref{3b5}) the bounds (\ref{App}) follow from 
\begin{equation} \label{dsy}
\sum_{n=K}^{N-1} \sum_{s=2}^{\infty} y_{ns} = \left\{ \begin{array}{ll} O(\Delta_K^{-4} K^{1-q} (\log N)^2) & \mbox{ if } K \rightarrow \infty, \cr
o(\log N) & \mbox{ if } K \mbox{ is fixed.} \end{array} \right.
\end{equation}

Proof of (\ref{dsy}) for $K \rightarrow \infty$:
Let $\lambda = \frac{2}{\log(1+\xi^2)}$.
We show that
\begin{eqnarray} \label{y1}
\sum_{n=K}^{N-1} y_{n2} & \leq & 2K^{1-q} (\log N)^2, \\ \label{y2}
\sum_{n=K}^{N-1} \sum_{3 \leq s \leq \lambda b_{nK}} y_{ns} & \leq & 2 \lambda^{\frac{1}{2}} K^{1-q} (\log N)^{\frac{3}{2}}, \\ \label{y3}
\sum_{n=K}^{N-1} \sum_{s > \lambda b_{nK}} y_{ns} & \leq & \Big[ \tfrac{(1+\xi^2)^{\frac{1}{2}}}{1-(1+\xi^2)^{-\frac{1}{2}}} \Big] K^{1-q} \log N.
\end{eqnarray}

Since $x^{-\frac{1}{2}} \log x$ is monotonically decreasing for $x \geq e^2$ and $b_{nK} \geq 1$ for $K$ large,
by (\ref{yn2}),
\begin{equation} \label{ineq1}
y_{n2} \leq \tfrac{\log[\exp(2b_{nK})]}{[\exp(2b_{nK})]^{\frac{1}{2}}} = \tfrac{2b_{nK}}{\exp(b_{nK})} = \tfrac{2K^{1-q}}{n} \log(\tfrac{n}{K^{1-q}}),
\end{equation}
and (\ref{y1}) holds. 
By (\ref{yns}),
$y_{ns} \leq (s-1)^{-\frac{1}{2}} \exp(-b_{nK})$,
and (\ref{y2}) follows from 
\begin{equation} \label{ineq2}
\sum_{3 \leq s \leq \lambda b_{nK}} y_{ns} \leq 2 (\lambda b_{nK})^{\frac{1}{2}} \exp(-b_{nK}) \leq 2 \lambda^{\frac{1}{2}} [\log(\tfrac{n}{K^{1-q}})]^{\frac{1}{2}} \tfrac{K^{1-q}}{n}.
\end{equation}
By (\ref{yns}),
$y_{ns} \leq (1+\xi^2)^{-\frac{s-1}{2}}$,
and (\ref{y3}) follows from
\begin{equation} \label{ineq3}
\sum_{s > \lambda b_{nK}} y_{ns} \leq \tfrac{(1+\xi^2)^{-\frac{\lambda b_{nK}-1}{2}}}{1-(1+\xi^2)^{-\frac{1}{2}}} = \Big[ \tfrac{(1+\xi^2)^{\frac{1}{2}}}{1-(1+\xi^2)^{-\frac{1}{2}}} \Big] \exp(-b_{nK}). 
\end{equation}

Proof of (\ref{dsy}) for $K$ fixed:
By the first two relations in (\ref{ineq1}),
$$y_{n2} \leq \tfrac{2 \{ \log(n/K)+\alpha \log[1+\log(n/K)] \}}{(n/K)[1+\log(n/K)]^{\alpha}} \leq 2(1+\alpha)[1+\log(\tfrac{n}{K})]^{1-\alpha} \tfrac{K}{n},
$$
and $\sum_{n=K}^{N-1} y_{n2} = o(\log N)$ follows from $\alpha>1$.
By the first inequality in (\ref{ineq2}),
\begin{eqnarray*}
\sum_{3 \leq s \leq \lambda b_{nK}} y_{ns} & \leq & \tfrac{2 \lambda^{\frac{1}{2}} \{ \log(n/K)+\alpha \log[1+\log(n/K)] \}^{\frac{1}{2}}}{(n/K)[1+\log(n/K)]^{\alpha}}, \cr
& \leq & 2 \lambda^{\frac{1}{2}}(1+\alpha)[1+\log(\tfrac{n}{K})]^{\frac{1}{2}-\alpha} \tfrac{K}{n},
\end{eqnarray*}
and $\sum_{n=K}^{N-1} \sum_{3 \leq s \leq \lambda b_{nK}} y_{ns} = o(\log N)$ follows from $\alpha > \frac{1}{2}$.
By (\ref{ineq3}),
$$\sum_{s > \lambda b_{nK}} y_{ns} \leq \Big[ \tfrac{(1+\xi^2)^{\frac{1}{2}}}{1-(1+\xi^2)^{-\frac{1}{2}}} \Big] [1+\log(\tfrac{n}{K})]^{-\alpha} \tfrac{K}{n},
$$
and $\sum_{n=K}^{N-1} y_{ns} = o(\log N)$ follows from $\alpha>0$.
\end{appendix}

\end{document}